\documentclass[a4paper,12pt]{amsart}
\usepackage{amssymb}
\usepackage{amsmath}
\usepackage{ascmac}
\usepackage{comment}
\usepackage{wrapfig}
\usepackage{txfonts}
\newtheorem{rem}{Remark}

\def\o#1{\overline{#1}}
\def\u#1{\underline{#1}}
\def\C{{\mathbb C}}

\def\Z{{\mathbb Z}}
\def\dfrac#1#2{{\displaystyle\frac{#1}{#2}}}
\def\sfrac#1#2{{#1}/{#2}}

\makeatletter
\@addtoreset{equation}{section}

\makeatother
\title{The Pad\'e interpolation method applied to $q$-Painlev\'e equations II (differential grid version)}
\author{Hidehito Nagao}
\address{Department of Arts and Science, National Institute of Technology, Akashi College, Hyogo 674-8501, Japan}
\email{nagao@akashi.ac.jp}
\keywords{Pad\'e method, Pad\'e interpolation, $q$-Painlev\'e equation.}
\subjclass[2010]{33D15, 34M55, 39A13, 41A21}

\begin{document}

\begin{abstract}
Recently we studied Pad\'e interpolation problems of $q$-grid, related to $q$-Painlev\'e equations of type $E_7^{(1)}$, $E_6^{(1)}$, $D_5^{(1)}$, $A_4^{(1)}$ and $(A_2+A_1)^{(1)}$. By solving those problems, we could derive evolution equations, scalar Lax pairs and determinant formulae of special solutions for the corresponding $q$-Painlev\'e equations. It is natural that the $q$-Painlev\'e equations were derived by the interpolation method of $q$-grid, but it may be interesting in terms of differential grid that the Pad\'e interpolation method of differential grid (i.e. Pad\'e approximation method) has been applied to the $q$-Painlev\'e equation of type $D_5^{(1)}$ by Y. Ikawa. In this paper we continue the above study and apply the Pad\'e approximation method to the $q$-Painlev\'e equations of type $E_6^{(1)}$, $D_5^{(1)}$, $A_4^{(1)}$ and $(A_2+A_1)^{(1)}$. Moreover determinant formulae of the special solutions for $q$-Painlev\'e equation of type $E_6^{(1)}$ are given in terms of the terminating $q$-Appell Lauricella function.
\end{abstract}

\maketitle
\renewcommand\baselinestretch{1.2}

\section{Introduction}\label{sec:intro}

\noindent
 In this paper we continue \cite{Nagao15} and apply the Pad\'e approximation method to
 the $q$-Painlev\'e equations of type $E_6^{(1)}$, $D_5^{(1)}$, $A_4^{(1)}$ and $(A_2+A_1)^{(1)}$.
 
\subsection{The background of discrete Painlev\'e equations}　\\

\noindent
In Sakai's theory \cite{Sakai01} the discrete Painlev\'e equations were classified on the basis of rational surfaces connected to extended affine Weyl groups. There exist three types of discrete Painlev\'e equations in the classification: elliptic difference ($e$-), multiplicative difference ($q$-) and additive difference ($d$-). The discrete Painlev\'e equations of $q$-difference type are classified as follows:

{\arraycolsep=1pt
\[
\begin{array}{ccccccccccccccccccc}
&&&&&&&&&&&&&&&{\mathbb Z}\\[-2mm]
&&&&&&&&&&&&&&\nearrow&&\searrow\\
&E_8^{(1)}
&\rightarrow&E_7^{(1)}&\rightarrow&E_6^{(1)}
&\rightarrow &\underset{(\mbox{$q$-$P_{\rm VI}$})}{D_5^{(1)}}&\rightarrow&\underset{(\mbox{$q$-$P_{\rm V}$})}{A_4^{(1)}}
&\rightarrow&\underset{(\mbox{$q$-$P_{\rm IV}$,$q$-$P_{\rm III}$})}{(A_2+A_1)^{(1)}}&\rightarrow&\underset{(\mbox{$q$-$P_{\rm II}$})}{(A_1+A_1^{\prime})^{(1)}}
&\rightarrow&\underset{(\mbox{$q$-$P_{\rm I}$})}{A_1^{(1)}}& &{\mathcal D}_6\\[5mm]
\\
\end{array}
\]
}
Here $A \to B$ means that $B$ is obtained from $A$ by degeneration.

\subsection{The background of the Pad\'e method}　\\

Pad\'e approximation/interpolation are closely related to Painlev\'e/Garnier equations. The Pad\'e method is a method for giving Painlev\'e equations, scalar Lax pairs and determinant formulae of special solutions simultaneously, by starting from suitable problems of Pad\'e approximation (of differential grid)/interpolation (of difference grid). In \cite{Yamada09-1} Y.Yamada has applied the Pad\'e method to continuous Painlev\'e equations of type $P_{\rm VI}$, $P_{\rm V}$, $P_{\rm IV}$ and Garnier system by using differential grid (i.e. Pad\'e approximation).
 
The Pad\'e method for discrete Painlev\'e equations has been applied to the following types: 
\begin{equation}\label{eq:pade-results}
\begin{tabular}{|c|c|c|c|c|c|c|c|c|}
\hline
&$e$-$E_8^{(1)}$&$q$-$E_8^{(1)}$&$q$-$E_7^{(1)}$&$q$-$E_6^{(1)}$&\multicolumn{2}{|c|}{$q$-$D_5^{(1)}$}&$q$-$A_4^{(1)}$&$q$-$(A_2+A_1)^{(1)}$\\
\hline\hline
&\cite{NTY13}&\cite{Yamada14}&\cite{Nagao15}&\cite{Ikawa13, Nagao15}&\cite{Ikawa13}&\cite{Nagao15}&\cite{Nagao15}&\cite{Nagao15}\\
\hline
grid &elliptic & $q$-quadric & $q$ & $q$ & differential & $q$ & $q$ & $q$\\
\hline
\end{tabular}
\end{equation} 
It is natural that the continuous/discrete Painlev\'e equations were derived by the interpolation of differential/difference grid respectively.
Here it may be interesting to note that the Pad\'e approximation method of differential grid (i.e. Pad\'e approximation) has been applied to the type $q$-$D_5^{(1)}$ in \cite{Ikawa13}. In this paper differential grid is applied to type $q$-$E_6^{(1)}$, $q$-$D_5^{(1)}$, $q$-$A_4^{(1)}$ and $q$-$(A_2+A_1)^{(1)}$.

\begin{rem} \label{rem:key}
{\rm {\bf On the key points of the Pad\'e method}

There are two key points to apply the Pad\'e approximation/interpolation method \cite{Ikawa13, Nagao15, NTY13, Yamada14}. The first key point is the appropriate choice of approximated/interpolated functions (see Table (\ref{eq:Ylistap}) and Remark \ref{rem:choice}). The second key point is to consider two linear $q$-difference three term relations (\ref{eq:L2L3matrixap}) satisfied by the error terms of the Pad\'e approximation/interpolation problems. Then the error terms can be expressed in terms of special solutions of $q$-Painlev\'e equations. Therefore the $q$-difference relations are the main subject in our study, and they naturally give the evolution equations, the Lax pairs and the special solutions for the corresponding $q$-Painlev\'e equations.$\square$
}
\end{rem}

\begin{rem} \label{rem:SOP}
{\rm {\bf On a connection between the Pad\'e method and the theory of semiclassical orthogonal polynomials}

The connection between semiclassical orthogonal polynomials (classical orthogonal polynomials related to a suitable weight function) and Painlev\'e/Garnier systems has been demonstrated in \cite{Magnus95}. It has been shown that coefficients of three term recurrence relations, satisfied by several semiclassical orthogonal polynomials, can be expressed in terms of solutions of Painlev\'e/Garnier systems (see \cite{Clarkson13, Nakazono13, OWF11,Van07,Witte09,Witte15,WO12} for example). Thus there exists a close connection between the Pad\'e method and the theory of semiclassical orthogonal polynomials. Namely, using both approaches, we can obtain the evolution equations, the Lax pairs and the special solutions for the corresponding Painlev\'e/Garnier systems. (The theory of semiclassical orthogonal polynomials is more general and the Pad\'e method is simpler. For example their relation was briefly proved in \cite{Yamada09-1}.)
$\square$
}
\end{rem}

\subsection{The purpose and the organization of this paper}\label{subsec:purpose}　\\

\noindent
The purpose of this paper is to apply the Pad\'e approximation method to type $q$-$E_6^{(1)}$, $q$-$D_5^{(1)}$, $q$-$A_4^{(1)}$ and $q$-$(A_2+A_1)^{(1)}$.  As the main results given in Section \ref{sec:approximation results}, the following items are presented for each type.

\vspace{3mm}
{\bf(a)} Setting of the Pad\'e approximation problem,

{\bf(b)} Contiguity relations,

{\bf(c)} The Painlev\'e equation,

{\bf(d)} The Lax pair,

{\bf(e)} Special solutions.\\

\noindent
This paper is organized as follows: In Section \ref{sec:approximation method} we explain the Pad\'e approximation method applied to the $q$-Painlev\'e equations, namely the methods for the items (a)--(e) above. In Section \ref{sec:approximation results} we present these main results for type $q$-$E_6^{(1)}$, $q$-$D_5^{(1)}$, $q$-$A_4^{(1)}$ and $q$-$(A_2+A_1)^{(1)}$. In Section \ref{sec:conc} we give a summary and discuss some future problems.
\section{Pad\'e approximation method of differential grid}\label{sec:approximation method}　

\noindent
In this section we explain the methods for deriving the items (a)--(e) in the main results given in Section \ref{sec:approximation results}. These contents of the items (b)--(d) below (i.e. Subsection \ref{subsec:itemb}, \ref{subsec:itemc} and \ref{subsec:itemd}) are almost the same as the items (b)--(d) in Section 2 of \cite{Nagao15}.

\subsection{(a) Setting of the Pad\'e approximation problem}\label{subsec:item(a)}　\\

\noindent
Let us consider the following approximation problem (of differential grid):

For a given function $Y(x)$, we look for functions $P_m(x)$ and $Q_n(x)$ which are polynomials of degree $m$ and $n$ $\in \Z_{\geq 0}$, satisfying the approximation condition 
\begin{equation}\label{eq:padeap}
Y(x)\equiv \frac{P_m(x)}{Q_n(x)}\quad(\mathrm{mod} \; x^{m+n+1}).
\end{equation}
We call this problem the "Pad\'e approximation problem (of differential grid)". Then the function $Y(x)$ is called the "generating function" (because $Y(x)$ generates the coefficients $p_k$ in power series (\ref{eq:padeYap}) in the item (e) below, i.e. Subsection \ref{subsec:iteme}), and the polynomials $P_m(x)$ and $Q_n(x)$ are called "approximating polynomials" respectively. The explicit expressions of the polynomials $P_m(x)$ and $Q_n(x)$, which are used in the computations for the item (e) above, are given in the formulae (\ref{eq:padePQap}) and (\ref{eq:padePQap2}) (see the item (e) below).

\begin{rem}\label{rem:normalizationap}
{\rm {\bf On the common normalization factor of the polynomials $P_m(x)$ and $Q_n(x)$ }

The common normalization factor of the approximating polynomials $P_m(x)$ and $Q_n(x)$ is not determined by the condition (\ref{eq:padeap}). However this normalization factor is not essential for our arguments, i.e. the main results in Section \ref{sec:approximation results} (see Remark \ref{rem:gaugeap}). $\square$
}
\end{rem}

\noindent
Fix a complex parameter $q$ ($0<|q|<1$). Let $a_i, b_i \in \C^{\times}$ be complex parameters.  
We establish the approximation problems (\ref{eq:padeap}) by specifying the generating functions $Y(x)$ as follows:
\begin{equation}\label{eq:Ylistap}
\begin{tabular}{|c||c|c|c|c|}
\hline
&$q$-$E_6^{(1)}$&$q$-$D_5^{(1)}$&$q$-$A_4^{(1)}$&$q$-$(A_2+A_1)^{(1)}$\\
\hline\hline
$\begin{array}{c}\\Y(x)\\ \\\end{array}$&
$\displaystyle\prod_{i=1}^3\dfrac{(a_ix;q)_{\infty}}{(b_ix;q)_{\infty}}$&
$\displaystyle\prod_{i=1}^2\dfrac{(a_ix;q)_{\infty}}{(b_ix;q)_{\infty}}$&
$\dfrac{(a_1x, a_2x;q)_{\infty}}{(b_1x;q)_{\infty}}$&
$(a_1x, a_2x;q)_{\infty}$\\
&$\dfrac{a_1a_2a_3q^m}{b_1b_2b_3q^n}=1$&&&\\
\hline
$q$-HGF&${}_3\varphi_2$&${}_2\varphi_1$&${}_2\varphi_1$&${}_1\varphi_1$\\
\hline
\end{tabular}
\end{equation}
Here $\dfrac{a_1a_2a_3q^m}{b_1b_2b_3q^n}=1$ is a constraint for the parameters in the case $q$-$E_6^{(1)}$, and the $q$-shifted factorials are defined by 
\begin{equation}\label{eq:qPochap}
\begin{array}{l}
(a_1, a_2, \cdots, a_i;q)_j:=\displaystyle\prod_{k=0}^{j-1}(1-a_1q^k)(1-a_2q^k)\cdots(1-a_iq^k)
\end{array}
\end{equation}
and the $q$-HGFs (i.e. the $q$-hypergeometric functions \cite{GaR04}) defined by 
\begin{equation}\label{eq:qHGFap}
\begin{array}{l}
{}_k\varphi_l\left(
\begin{array}{ccc}
a_1,&\cdots,&a_k\\[0mm]
b_1,&\cdots,&b_l
\end{array}
; q, x
\right)
:=\displaystyle\sum_{s=0}^{\infty}\dfrac{(a_1,\cdots,a_k ; q)_s}{(b_1,\cdots,b_l,q ; q)_s}\left[(-1)^sq^{\left(\substack{s\\2}\right)}\right]^{1+l-k}x^s
\end{array}
\end{equation}
with $\left(\substack{s\\2}\right)=\sfrac{s(s-1)}{2}$.

The hypergeometric solutions to the $q$-Painlev\'e equations of type $E_6^{(1)}$, $D_5^{(1)}$, $A_4^{(1)}$ and $(A_2+A_1)^{(1)}$ were given in terms of the $q$-hypergeometric functions ${}_3\varphi_2$, ${}_2\varphi_1$, ${}_2\varphi_1$, and ${}_1\varphi_1$ respectively in \cite{KMNOY04}. The functions $Y(x)$ in Table (\ref{eq:Ylistap}) generate the coefficients $p_k$ in power series (\ref{eq:padeYap}) in terms of the terminating $q$-hypergeometric functions ${}_3\varphi_2$ (\ref{eq:E6pap}) (see Remark \ref{rem:trans}), ${}_2\varphi_1$ (\ref{eq:D5pap}), ${}_2\varphi_1$ (\ref{eq:A4pap}), and ${}_1\varphi_1$ (\ref{eq:A21pap}). 

\begin{rem} \label{rem:choice}
{\rm {\bf On the choice of the generating functions $Y(x)$}

\noindent
One may wonder how the generating functions $Y(x)$ are appropriately chosen. However there is no theoretical choice of the functions $Y(x)$ in the Pad\'e approximation method as far as we know. We chose the functions $Y(x)$ to type $E_6^{(1)}$, $A_4^{(1)}$ and $(A_2+A_1)^{(1)}$ by a extension and reductions from the generating function to type $D_5^{(1)}$ given in \cite{Ikawa13}. 
$\square$
}
\end{rem}

Let us consider yet another Pad\'e problem where some parameters $a_i,b_i, m$ and $n$ in the generating functions $Y(x)$ are shifted. The parameter shift operators $T$ are given as follows:
\begin{equation}\label{eq:Tlistap}
\begin{tabular}{|c||ccc|}
\hline
&parameter&&\\
\hline\hline
$q$-$E_6^{(1)}$&$(a_1,a_2,a_3,b_1,b_2,b_3,m,n)$&$\mapsto$&$(qa_1,a_2,a_3,qb_1,b_2,b_3,m,n)$\\
\hline
$q$-$D_5^{(1)}$&$(a_1,a_2,b_1,b_2,m,n)$&$\mapsto$&$(qa_1,a_2,qb_1,b_2,m,n)$\\
\hline
$q$-$A_4^{(1)}$&$(a_1,a_2,b_1,m,n)$&$\mapsto$&$(qa_1,a_2,qb_1,m,n)$\\
\hline
$q$-$(A_2+A_1)^{(1)}$&$(a_1,a_2,m,n)$&$\mapsto$&$(qa_1,a_2,m,n)$\\
\hline
\end{tabular}
\end{equation}
Here the operators $T$ are called the "time evolutions", because they specify the directions of the time evolutions for $q$-Painlev\'e equations.\\

\subsection{(b) Contiguity relations}\label{subsec:itemb}　\\

\noindent
Let us consider two linear three term relations: $L_2(x)=0$ between $y(x), y(qx), \o{y}(x)$ and $L_3(x)=0$ between $y(x), \o{y}(x), \o{y}(x/q)$ satisfied by fundamental solutions $y(x)=P_m(x)$, $Y(x)Q_n(x)$, where $L_2$ and $L_3$ are given as expressions 
\begin{align}\label{eq:L2L3matrixap}
L_2(x)\propto
\begin{vmatrix}\nonumber 
y(x) & y(qx) & \o{y}(x) \\
P_m(x) & P_m(qx) & \o{P}_m(x)\\
Y(x)Q_n(x) & Y(qx)Q_n(qx) & \o{Y}(x)\o{Q}_n(x)
\end{vmatrix}, \\
L_3(x)\propto
\begin{vmatrix} 
y(x) & \o{y}(x) & \o{y}(x/q) \\
P_m(x) & \o{P}_m(x) & \o{P}_m(x/q)\\
Y(x)Q_n(x) & \o{Y}(x)\o{Q}_n(x) & \o{Y}(x/q)\o{Q}_n(x/q)
\end{vmatrix}.
\end{align}
Then the linear relations $L_2=0$ and $L_3=0$ are called the "contiguity relations", and the contiguity relations are the main subject in our study. Here for any object $F$ the corresponding shifts are denoted as $\o{F}:=T(F)$ and $\u{F}:=T^{-1}(F)$, and the shift operator $T$ acts on parameters given in Table (\ref{eq:Tlistap}).\\

%

\noindent
We show the method of computation of the contiguity relations $L_2=0$ and $L_3=0$.\\
Set ${\bf y}(x):=\left[\begin{array}{c}P_m(x)\\Y(x)Q_n(x)\end{array}\right]$ and define Casorati determinants $D_i(x)$ by
\begin{equation}\label{eq:Casoratiap}
\begin{array}{l}
D_1(x):=\det[{\bf y}(x),{\bf y}(qx)],\hspace{3mm} D_2(x):=\det[{\bf y}(x),{\o{\bf y}}(x)],\hspace{3mm} D_3(x):=\det[{\bf y}(qx),\o{{\bf y}}(x)].
\end{array}
\end{equation}
Then the expressions (\ref{eq:L2L3matrixap}) can be rewritten as follows:
\begin{equation}\label{eq:L2L3ap}
\begin{array}{l}
L_2(x)\propto D_1(x) \o{y}(x)-D_2(x)y(qx)+D_3(x)y(x),
\\
L_3(x)\propto \o{D}_1(\sfrac{x}{q})y(x)+D_3(\sfrac{x}{q}) \o{y}(x)-D_2(x)\o{y}(\sfrac{x}{q}).
\end{array}
\end{equation}
Define basic quantities $G(x), K(x)$ and $H(x)$ (e.g. (\ref{eq:E6GKHap}), (\ref{eq:D5GKHap})) by
\begin{equation}\label{eq:GKHap}
\begin{array}{l}
G(x):=\sfrac{Y(qx)}{Y(x)},\quad K(x):=\sfrac{\o{Y}(x)}{Y(x)},\quad H(x):={\rm L.C.M}(G_{\rm den}(x), K_{\rm den}(x)).
\end{array}
\end{equation}
Here $G_{\rm den}(x)$ and $G_{\rm num}(x)$ are defined as the polynomials of the denominator and the numerator in $G(x)$ respectively, and $K_{\rm den}(x)$ and $K_{\rm num}(x)$ are similarly defined. For example in case of $q$-$E_6^{(1)}$, $G_{\rm den}(x)=\prod_{i=1}^3 (1-a_ix)$, $G_{\rm num}(x)=\prod_{i=3}^3 (1-b_ix)$, $K_{\rm den}(x)=1-a_1x$, $K_{\rm num}(x)=1-b_1x$ (see eq.(\ref{eq:E6GKHap})). Substituting these quantities into the determinants (\ref{eq:Casoratiap}), we obtain the following expressions:
\begin{equation}\label{eq:Drelationap}
\begin{array}{l}
D_1(x)=\dfrac{Y(x)}{G_{\rm den}(x)}\left\{G_{\rm num}(x)P_{m}(x)Q_n(qx)-G_{\rm den}(x)P_{m}(qx)Q_n(x) \right\},
\\[5mm]
D_2(x)=\dfrac{Y(x)}{K_{\rm den}(x)}\left\{K_{\rm num}(x)P_{m}(x)\o{Q}_n(x)-K_{\rm den}(x)\o{P}_{m}(x)Q_n(x) \right\},
\\[5mm]
D_3(x)=\dfrac{Y(x)}{H(x)}
\left\{\dfrac{H(x)}{K_{\rm den}(x)}K_{\rm num}(x)P_{m}(qx)\o{Q}_n(x)-\dfrac{H(x)}{G_{\rm den}(x)}G_{\rm num}(x)\o{P}_{m}(x)Q_n(qx)\right\}.
\end{array}
\end{equation}
Using the approximation condition (\ref{eq:padeap}) and the form of the basic quantities $G(x)$, $K(x)$, $H(x)$ (e.g. eqs.(\ref{eq:E6GKHap}), (\ref{eq:D5GKHap})), we can investigate positions of zeros (e.g. $x=0$) and degrees of the polynomials (e.g. $G_{\rm num}(x)P_{m}(x)Q_n(qx)-G_{\rm den}(x)P_{m}(qx)Q_n(x)$) within braces $\{\hspace{3mm}\}$ of the expressions (\ref{eq:Drelationap}). Then we can simply compute the determinants $D_i(x)$ (e.g. eqs.(\ref{eq:E6Dap}), (\ref{eq:D5Dap}))  except for some factors such as $1-fx$, $1-\sfrac{x}{g}$ and $c_i$ in $D_i(x)$, where $f$, $g$ and $c_i$ are constants with respect to $x$ (see Remark \ref{rem:genericap}). In this way we obtain the contiguity relations $L_2=0$ and $L_3=0$ (e.g. eqs.(\ref{eq:E6L2L3ap}), (\ref{eq:D5L2L3ap})).

\begin{rem}\label{rem:gaugeap}
{\rm {\bf On the gauge invariance of the product $C_0C_1$}

\noindent
When the common normalization factor of the approximating polynomials $P_m(x)$ and $Q_n(x)$ is changed, an $x$-independent gauge transformation of $y(x)$ is induced in the contiguity relations $L_2=0$ and $L_3=0$. Under the $x$-independent gauge transformation of $y(x)$: $y(x)\mapsto Gy(x)$, the coefficients of $\o{y}(x), y(x/q), y(x)$ and $y(x), \o{y}(x), \o{y}(x/q)$ in $L_2=0$ and $L_3=0$ (\ref{eq:L2L3ap}) change as follows:
\begin{equation}
\begin{array}{l}
(D_1(x) : D_2(x) : D_3(x))\mapsto(\sfrac{\o{G}D_1(x)}{G} : D_2(x) : D_3(x))\\
(\o{D}_1(x/q) : D_3(x/q) : D_2(x))\mapsto(\sfrac{G\o{D}_1(x/q)}{\o{G}} : D_3(x/q) : D_2(x)).
\end{array}
\end{equation}
The coefficients $C_0$ and $C_1$ in $L_2=0$ and $L_3=0$ (e.g. eqs.(\ref{eq:E6L2L3ap}), (\ref{eq:D5L2L3ap})) are defined as the normalization factors of the coefficients of $\o{y}(x)$ and $y(x)$ respectively. Then $C_0$ and $C_1$ change under the gauge transformation, but the product $C_0C_1$ is a gauge invariant quantity. Moreover $C_0$ and $C_1$ do not appear in the final form of the $q$-Painlev\'e equations (e.g. eqs.(\ref{eq:E6eqap}), (\ref{eq:D5eqap})). $\square$
}
\end{rem}

\begin{rem}\label{rem:genericap}
{\rm {\bf On two meanings of variables $f, g$ and parameters $m, n$}

\noindent
We use $f$ and $g$ with two different meanings. The first meaning is constants (i.e. special solutions) $f$ and $g$ which are explicitly determined in terms of parameters $a_i, b_i, m$ and $n$ by the Pad\'e approximation problem (e.g. eqs.(\ref{eq:E6Dap}), (\ref{eq:E6L2L3ap}), (\ref{eq:E6solap}), (\ref{eq:D5Dap}), (\ref{eq:D5L2L3ap}), (\ref{eq:D5solap})). The second meaning is generic variables (i.e. generic solutions) $f$ and $g$ apart from the Pad\'e approximation problem (e.g. eqs.(\ref{eq:E6eqap}), (\ref{eq:E6L1L2ap}), (\ref{eq:D5eqap}), (\ref{eq:D5L1L2ap})), namely $f$ and $g$ are unknown functions in the $q$-Painlev\'e equation. In the items (c), (d) (resp. in the items (b), (e)) we consider $f$ and $g$ in the second meaning (resp. in the first meaning).

Similarly we use $m$ and $n$ with two meanings. In the first meaning $m$ and $n$ $\in \Z_{\geq 0}$ are integer parameters  (e.g. eqs.(\ref{eq:E6Yap}), (\ref{eq:E6Dap}), (\ref{eq:E6L2L3ap}), (\ref{eq:D5Yap}), (\ref{eq:D5Dap}), (\ref{eq:D5L2L3ap})). In the second meaning $m$ and $n$ $\in \C^{\times}$ are generic complex parameters, namely $q^m$ and $q^n$ are replaced by generic parameters $a_0$ and $b_0$ respectively 
(e.g. eqs.(\ref{eq:E6eqap}), (\ref{eq:E6L1L2ap}), (\ref{eq:D5eqap}), (\ref{eq:D5L1L2ap})). In the items (c), (d) (resp. in the items (a), (b), (e)) we consider $f$ and $g$ in the second meaning (resp. in the first meaning). Then the result of the compatibility of the contiguity relations $L_2=0$ and $L_3=0$ also holds with respect to the second meaning. $\square$
}
\end{rem}

\subsection{(c) The $q$-Painlev\'e equation}\label{subsec:itemc}　\\

\noindent
Let us consider generic variables $f, g$ and generic parameter $a_0, b_0$ as in the second meaning in Remark \ref{rem:genericap},  we can derive the $q$-Painlev\'e equation as the necessary condition for the compatibility of the contiguity relations $L_2=0$ and $L_3=0$ (e.g. eqs.(\ref{eq:E6L2L3ap}), (\ref{eq:D5L2L3ap})). Computing the compatibility condition, we determine three variables $\u{g}, \o{f}$ and $C_0C_1$. Expressions for variables $\u{g}$ and $\o{f}$ are obtained in terms of variables $f$ and $g$. An expression for the product $C_0C_1$ is obtained in terms of variables $f, g$ and $\o{f}$ (and hence in terms of variables $f$ and $g$).

The first and the second expressions are the $q$-Painlev\'e equation (e.g. eqs.(\ref{eq:E6eqap}), (\ref{eq:D5eqap})). The third expression is a constraint for the product $C_0C_1$ (e.g. eqs.(\ref{eq:E6C0C1ap}), (\ref{eq:D5C0C1ap})).

\subsection{(d) The Lax pair}\label{subsec:itemd}　\\

\noindent
Let us consider two linear three term equations for the unknown function $y(x)$: $L_1(x)=0$ between $y(qx), y(x), y(x/q)$ and $L_2(x)=0$ between $y(x), y(qx), \o{y}(x)$, where $L_1$ and $L_2$ are given as expressions 
\begin{equation}\label{eq:L1L2ap}
\begin{array}{l}
L_1(x)=A_1(x)y(\sfrac{x}{q})+A_2(x)y(x)+A_3(x)y(qx),
\\
L_2(x)=A_4(x)\o{y}(x)+A_5(x)y(x)+A_6(x)y(qx).
\end{array}
\end{equation}  
The linear three term equations $L_1=0$ and $L_2=0$ (\ref{eq:L1L2ap}) are called the "scalar Lax pair", 
 when the compatibility condition of the linear equations $L_1=0$ and $L_2=0$ (\ref{eq:L1L2ap}) is equivalent to a $q$-Painlev\'e equation.

We present how to compute the scalar Lax pair. Similarly to the item (c), let us consider generic variables $f, g$ and generic parameter $a_0, b_0$ as in the second meaning in Remark \ref{rem:genericap}. The Lax pair $L_1=0$ and $L_2=0$, which satisfies the compatibility condition, is derived using the results of the items (a)--(c) as follows: The Lax equation $L_2=0$ (e.g. eqs.(\ref{eq:E6L1L2ap}), (\ref{eq:D5L1L2ap})) in the item (d) is the same as the contiguity relation $L_2=0$ (e.g. eqs.(\ref{eq:E6L2L3ap}), (\ref{eq:D5L2L3ap})) in the item (b) under an $x$-independent gauge transform of $y(x)$ and changes of parameters. We can obtain the Lax equation $L_1=0$ as follows: First combining the contiguity relations $L_2=0$ and $L_3=0$ (e.g. eqs.(\ref{eq:E6L2L3ap}), (\ref{eq:D5L2L3ap})) under generic variables $f, g$ and generic parameters $a_0, b_0$, one obtains a linear equation between the three terms $y(q x), y(x)$ and $y(x/q)$ (See the figure below), whose coefficient functions depend on the variables $f, g, \o{f}, C_0$ and $C_1$. However the variables $C_0$ and $C_1$ appear through the product $C_0C_1$. Therefore expressing $\o{f}$ (e.g. eqs.(\ref{eq:E6eqap}), (\ref{eq:D5eqap})) and $C_0C_1$ (e.g. eqs.(\ref{eq:E6C0C1ap}), (\ref{eq:D5C0C1ap})) in terms of $f$ and $g$ only, one obtains the Lax equation $L_1=0$ (e.g. eqs.(\ref{eq:E6L1L2ap}), (\ref{eq:D5L1L2ap})). 

\begin{center}\setlength{\unitlength}{1.2mm}
\begin{picture}(50,28)(-5,-3)
\put(-1,23){$\o{y}(\sfrac{x}{q})$}
\put(20,23){$\o{y}(x)$}
\put(-1,-4){$y(\sfrac{x}{q})$}
\put(20,-4){$y(x)$}
\put(41,-4){$y(qx)$}
\put(10,-4){$L_1(x)$}
\put(0,-1){\line(1,0){42}}
\put(0,0){\line(1,0){20}}
\put(0,0){\line(0,1){20}}
\put(0,20){\line(1,-1){20}}
\put(1,21){\line(1,0){20}}
\put(21,1){\line(0,1){20}}
\put(1,21){\line(1,-1){20}}
\put(2,4){$L_2(\sfrac{x}{q})$}
\put(11,14){$L_3(x)$}
\put(22,0){\line(1,0){20}}
\put(22,0){\line(0,1){20}}
\put(22,20){\line(1,-1){20}}
\put(26,4){$L_2(x)$}
\end{picture}
\end{center}

Then, under generic variables $f, g$ and generic parameters $a_0, b_0$, the $q$-Painlev\'e equation (e.g. eqs.(\ref{eq:E6eqap}), (\ref{eq:D5eqap})) is sufficient for the compatibility of the Lax pair $L_1=0$ and $L_2=0$ (e.g. eqs.(\ref{eq:E6L1L2ap}), (\ref{eq:D5L1L2ap})). The similar proofs have already been given in \cite{KNY15, Yamada09-2, Yamada11}.

\subsection{(e) Special solutions}\label{subsec:iteme}　\\

\noindent
By construction, expressions for constants $f$ and $g$ as in the first meaning in Remark \ref{rem:genericap} give a special solution for the $q$-Painlev\'e equation. We present how to compute determinant formulae of the special solutions $f$ and $g$.\\

\noindent
We derive the formulae (\ref{eq:padePQap}) which are convenient for computing the special solutions $f$ and $g$. We can assume the power series
\begin{equation}\label{eq:padeYap}
Y(x)=\displaystyle\sum_{k=0}^{\infty}p_kx^k,\quad p_0=1,\quad p_i=0 \quad(i<0)
\end{equation}
since the generating functions $Y(x)$ in Table (\ref{eq:Ylistap}) are holomorphic near $x=0$. Then for each type of the functions $Y(x)$, the coefficients $p_k$ (e.g. eqs.(\ref{eq:D5pap}), (\ref{eq:A4pap})) are given respectively as the terminating cases of $q$-hypergeometric functions defined by eq.(\ref{eq:qHGFap}).

For a given function $Y(x)$, the polynomials $P_m(x)$ and $Q_n(x)$ of degree $m$ and $n$ for the approximation condition (\ref{eq:padeap}) are given by the following determinant expressions:
\begin{equation}\label{eq:padePQap}
P_m(x)=\displaystyle\sum_{i=0}^{m}s_{(m^n,i)}x^i,\quad Q_n(x)=\displaystyle\sum_{i=0}^{n}s_{((m+1)^i,m^{n-i})}(-x)^i
\end{equation}
where $s_{\lambda}$ is the Schur function defined by  the Jacobi-Trudi formula
\begin{equation}\label{eq:Jacobi-Trudi}
s_{(\lambda_1,\cdots,\lambda_l)}:=\det(p_{\lambda_i-i+j})_{i,j=1}^l.
\end{equation}

We show the derivation of the expressions (\ref{eq:padePQap}) as follows: The approximating polynomial $Q_n(x)$ satisfying the condition (\ref{eq:padeap}) can be given  as the second expression of eq.(\ref{eq:padePQap}) by Cramer's rule. The approximating polynomial $P_m(x)$ satisfying the condition (\ref{eq:padeap}) is given by the following computation: By using the relation
\begin{equation}\label{eq:Relation}
x^n Y(x)=\displaystyle\sum_{k=0}^{\infty}p_kx^{k+n}=\displaystyle\sum_{k=0}^{\infty}p_{k-n}x^k,
\end{equation}
we have
\begin{align}
Y(x)Q_n(x)
=&
\begin{vmatrix}\label{eq:YQPmatrix1}
p_m  & p_{m+1}  & \cdots & p_{m+n} \\
\vdots & \ddots   & \ddots & \vdots \\
p_{m-n+1} & \cdots & p_m & p_{m+1}\\
x^nY(x)& \cdots & xY(x) & Y(x)
\end{vmatrix}
=&\displaystyle\sum_{k=0}^{\infty}
\begin{vmatrix} 
p_m  & p_{m+1}  & \cdots & p_{m+n} \\
\vdots & \ddots   & \ddots & \vdots \\
p_{m-n+1} & \cdots & p_m & p_{m+1}\\
p_{k-n} & \cdots & p_{k-1} & p_k
\end{vmatrix}
x^k
=\displaystyle\sum_{k=0}^{\infty}s_{(m^n,k)}x^k.
\end{align}
Here we note that
\begin{align}\label{eq:YQPmatrix2}
\displaystyle\sum_{k=m+1}^{m+n}s_{(m^n,k)}x^k=
0.
\end{align}
Substituting the relation (\ref{eq:YQPmatrix2}) into the expression (\ref{eq:YQPmatrix1}), we obtain
\begin{equation}\label{eq:YQPmatrix}
Y(x)Q_n(x)=\Big(\displaystyle\sum_{k=0}^{m}+\displaystyle\sum_{k=m+n+1}^{\infty}\Big)s_{(m^n,k)}x^k.
\end{equation}
Hence the desired polynomial $P_m(x)$ is given  as the first expression of the formulae (\ref{eq:padePQap}).

Furthermore the polynomials $P_m(x)$ and $Q_n(x)$ in the formulae (\ref{eq:padePQap}) can be expressed in terms of a single determinant as 
\begin{equation}\label{eq:padePQap2}
P_m (x)=x^ms_{(m^{n+1})}|_{p_i \rightarrow \sum_{j=0}^{i}x^{-j}p_{i-j}}, \quad Q_n (x)=(-x)^ns_{((m+1)^n)}|_{p_i \rightarrow p_i -x^{-1}p_{i-1}}.
\end{equation}

The formulae (\ref{eq:padePQap}) and (\ref{eq:padePQap2}) have already appeared in \cite{Yamada09-1}. 

Then we apply the general results described above to the case $N=1$ (\ref{eq:D5Yap}) and $N=2$ (\ref{eq:E6Yap}) of the function $\psi(x):=\prod_{i=1}^{N+1}\dfrac{(a_i x)_\infty}{(b_i x)_\infty}$, 
which can be written as
\begin{equation}\label{eq:Tsuda}
\psi(x)=\exp\Big(\sum_{k=1}^{\infty}\sum_{s=1}^{N+1}\dfrac{b_s^{k}-a_s^{k}}{k(1-q^k)}x^k\Big)=\sum_{k=0}^{\infty}p_k x^k.
\end{equation}
We note that this kind of expression (\ref{eq:Tsuda}) has already appeared in \cite{Tsuda10}.\\

\noindent
We show the method of computation of the special solutions $f$ and $g$.\\
The expressions for the special solutions $f$ and $g$ can be derived by comparing the determinants $D_i(x)$ in eq.(\ref{eq:Drelationap}) and $D_i(x)$ (e.g. eqs.(\ref{eq:E6Dap}), (\ref{eq:D5Dap})) in the item (b) as the identity with respect to the variable $x$ and applying the formulae (\ref{eq:padePQap}) and (\ref{eq:padePQap2}). 

For example the computation for the case $q$-$E_6^{(1)}$ is as follows: Substituting $x=\sfrac{1}{a_i}$ $(i=1, 2)$ into the determinants $D_1(x)$ in eq.(\ref{eq:Drelationap}) and $D_1(x)$ in (\ref{eq:E6Dap}) respectively, we obtain an expression for the special solution $f$ in the first equation of eq.(\ref{eq:E6solap}) by comparing the two expressions for $D_1(x)$ and applying the formulae (\ref{eq:padePQap2}). Similarly substituting $x=\sfrac{1}{b_i}$ $(i=2, 3)$ into the determinants $D_3(x)$ in eq.(\ref{eq:Drelationap}) and $D_3(x)$ in eq.(\ref{eq:E6Dap}) respectively, we obtain an expression for the special solution $g$ in the second equation of eq.(\ref{eq:E6solap}) by comparing the two expressions for $D_3(x)$ and applying the formulae (\ref{eq:padePQap2}).\\

\section{Main results}\label{sec:approximation results}　\\

\noindent
In this section for each case $q$-$E_6^{(1)}$,  $q$-$D_5^{(1)}$, $q$-$A_4^{(1)}$ and $q$-$(A_2+A_1)^{(1)}$, we present the results obtained through the method, which was explained in Section \ref{sec:approximation method}.

We use the following notations:
\begin{equation}\label{eq:notationsap}
\begin{array}{l}
\sfrac{a_1a_2\cdots a_n}{b_1b_2\cdots b_n}:=\dfrac{a_1a_2\cdots a_n}{b_1b_2\cdots b_n},\\[5mm]
\tau_{m,n}:=s_{(m^n)},\\[5mm]
T_{a_i}(F):=F |_{a_i \to qa_i},\quad T_{a_i}^{-1}(F):= F |_{a_i \to \sfrac{a_i}{q}}
\end{array}
\end{equation}
for any quantity (or function) $F$ depending on variables $a_i$ and $b_i$, and by definition (\ref{eq:Jacobi-Trudi}), the Schur function $s_{(m^n)}$ is expressed as  
\begin{equation}
s_{(m^n)}=
\begin{vmatrix} 
p_m  & p_{m+1} &\cdots & p_{m+n-1} \\
p_{m-1} & \ddots &  & \vdots \\
\vdots & \cdots & \ddots &\vdots \\
p_{m-n+1} & \cdots  & \cdots  & p_m
\end{vmatrix}
\end{equation}
where the element $p_k$ is defined in the power series (\ref{eq:padeYap}).\\

\subsection{Case $q$-$E_6^{(1)}$}\label{subsubsec:E6ap}　\\

\noindent
{\bf(a)} Setting of the Pad\'e approximation problem\\

In Table (\ref{eq:Ylistap}) the generating function and the constraint are established as
\begin{equation}\label{eq:E6Yap}
Y(x):=\dfrac{(a_1 x, a_2 x, a_3x;q)_\infty}{(b_1 x, b_2 x, b_3x;q)_\infty},\quad \dfrac{a_1a_2a_3q^m}{b_1b_2b_3q^n}=1, 
\end{equation}
and in Table (\ref{eq:Tlistap}) the time evolution is chosen as
\begin{equation}\label{eq:E6Tap}
T: (a_1,a_2,a_3,b_1,b_2,b_3,m,n) \mapsto (qa_1,a_2,a_3,qb_1,b_2,b_3,m,n).
\end{equation}
\\
\noindent
{\bf(b)} Contiguity relations\\

By the definition (\ref{eq:GKHap}) we have the basic quantities
\begin{equation}\label{eq:E6GKHap}
G(x)=\displaystyle\prod_{i=1}^{3}\dfrac{(1-b_ix)}{(1-a_ix)},\quad 
K(x)=\dfrac{1-b_1x}{1-a_1x},\quad H(x)=\prod_{i=1}^{3}(1-a_ix),
\end{equation}
and by the expression (\ref{eq:Drelationap}) we obtain the Casorati determinants
\begin{equation}\label{eq:E6Dap}
\begin{array}{l}
D_1(x)=:\sfrac{c_0(1-xf)x^{m+n+1}Y(x)}{\prod_{i=1}^{3}(1-a_ix)},\quad
D_2(x)=:\sfrac{c_1x^{m+n+1}Y(x)}{(1-a_1x)},
\\
D_3(x)=:\sfrac{c_1\dfrac{a_2a_3q^mg}{b_1}(1-b_1x)(1-\sfrac{x}{g})x^{m+n+1}Y(x)}{\prod_{i=1}^{3}(1-a_ix)}
\end{array}
\end{equation}
where $f, g, c_0$ and $c_1$ are constants depending on parameters $a_i, b_i \in \C^{\times}(i=1,2,3)$, $m,n \in \Z_{\geq 0}$ but independent of $x$. Then the contiguity relations $L_2=0$ and $L_3=0$ are expressed by
\begin{equation}\label{eq:E6L2L3ap}
\begin{array}{l}
L_2(x)=C_0(1-xf)\o{y}(x)-(1-a_2x)(1-a_3x)y(qx)
+\dfrac{a_2a_3q^mg}{b_1}(1-b_1x)(1-\sfrac{x}{g})y(x),\\
L_3(x)=C_1(1-\sfrac{x\o{f}}{q})y(x)+\dfrac{a_2a_3q^mg}{b_1}(1-a_1x)(1-\sfrac{x}{qg})\o{y}(x)
-q^{m+n+1}(1-\sfrac{b_2x}{q})(1-\sfrac{b_3x}{q})\o{y}(\sfrac{x}{q})
\end{array}
\end{equation}
where $C_0=\sfrac{c_0}{c_1}$ and $C_1=\sfrac{\o{c}_0}{c_1}$.\\

Take note that in the items (c) and (d) below we study the contiguity relations $L_2=0$ and $L_3=0$ (\ref{eq:E6L2L3ap}) for generic complex parameters $a_0=q^m$, $b_0=q^n$ $(m, n \in \C^{\times})$ and generic variables $f, g$ (depending on parameters $a_i, b_i \in \C^{\times}, i=0,1,2,3$) apart from the Pad\'e approximation problem (\ref{eq:padeap}) with eqs.(\ref{eq:E6Yap}) and (\ref{eq:E6Tap}). (see Remark \ref{rem:genericap})\\

\noindent
{\bf(c)} The $q$-Painlev\'e equation\\

Compatibility of the contiguity relations $L_2=0$ and $L_3=0$ (\ref{eq:E6L2L3ap}) gives the evolution equations and the constraint on the product $C_0C_1$ as follows:
\begin{equation}\label{eq:E6eqap}
\begin{array}{l}
(fg-1)(f\u{g}-1)=\dfrac{b_0b_1^2}{a_0a_2^2a_3^2}\dfrac{(f-a_2)(f-a_3)(f-b_2)(f-b_3)}{(f-a_1)(f-b_1)},
\\
(fg-1)(\o{f}g-1)=qa_1b_1\dfrac{(g-\sfrac{1}{a_2})(g-\sfrac{1}{a_3})(g-\sfrac{1}{b_2})(g-\sfrac{1}{b_3})}{(g-\sfrac{b_1}{a_0a_2a_3})(g-\sfrac{qb_0b_1}{a_2a_3})}
\end{array}
\end{equation}
and
\begin{equation}\label{eq:E6C0C1ap}
C_0C_1=\sfrac{a_0(b_1 - a_0a_2 a_3 g) (qb_0b_1-a_2 a_3 g)}{b_1^2}.
\end{equation}
The evolution equations (\ref{eq:E6eqap}) are equivalent to the $q$-Painlev\'e equation of type $E_6^{(1)}$ given in \cite{Ikawa13, KMNOY04, KNY15, Nagao15, RGTT01}.
The 8 singular points in coordinates $(f,g)$ are on the two lines $f=\infty$ and $g=\infty$ and one curve $fg=1$ as follows:
\begin{equation}
\begin{array}{l}
(f,g)=(\sfrac{1}{a_2}, a_2),(\sfrac{1}{a_3}, a_3),(\sfrac{1}{b_2}, b_2),(\sfrac{1}{b_3}, b_3),\\
\phantom{(f,g)=}(a_1, \infty),(b_1, \infty),(\infty,\sfrac{b_1}{a_0a_2a_3}),(\infty,\sfrac{qb_0b_1}{a_2a_3}).
\end{array}
\end{equation}

\noindent
{\bf(d)} The Lax pair\\

The contiguity relations $L_2=0$ and $L_3=0$ (\ref{eq:E6L2L3ap}) give two scalar Lax equations $L_1=0$ and $L_2=0$ expressed by
\begin{equation}\label{eq:E6L1L2ap}
\begin{array}{l}
L_1(x)=\dfrac{a_0b_0(1-b_1 x/q) (1-b_2 x/q)
   (1-b_3 x/q)}{(1-f
   x/q)}
   \Big[y\left(\sfrac{x}{q}\right)-\dfrac{a_1 (1-a_2 x/q)
   (1-a_3 x/q)}{b_0b_2 b_3
   (1-b_1 x/q) (g-x/q)} y(x)\Big]
   \\[5mm]\phantom{L_1(x):}
   +\dfrac{(1-a_1 x) (1-a_2 x)
   (1-a_3 x)}{q(1-f
   x)} \Big[y(q
   x)-\dfrac{b_0b_2 b_3(1-b_1
   x) (g-x)}{a_1 (1-a_2
   x) (1-a_3 x)} y(x)\Big]
   \\[5mm]\phantom{L_1(x):}
   +\dfrac{a_0a_1}{b_2 b_3 g}
   \Big[\left(1-\dfrac{b_2 b_3 g
   }{qa_0a_1}\right)
   \left(1-\dfrac{b_0b_2 b_3 g}{a_1}\right)
   +\dfrac{x
   (1-a_2 g) (1-a_3 g)
   (1-b_2 g) (1-b_3 g)}{(1-f g)
   (g q-x)}\Big]y(x),
   \\[5mm]
L_2(x)=(1-xf)\o{y}(x)-(1-a_2x)(1-a_3x)y(qx)
+\dfrac{a_0a_2a_3g}{b_1}(1-b_1x)(1-\sfrac{x}{g})y(x).
\end{array}
\end{equation}
The scalar Lax pair $L_1=0$ and $L_2=0$ (\ref{eq:E6L1L2ap}) is expected to be equivalent to the 2 $\times$ 2 matrix ones in \cite{Sakai06, WO12} and the scalar ones in \cite{Ikawa13, Nagao15, Yamada11} by using suitable gauge transformations of $y(x)$. (Note that there are some typographical errors in eqs.(30) and (31) in \cite{Ikawa13}, namely the expressions $(b_4q^{\prime}-z)$ and $Y(x)$ should read $(b_4q^{\prime}-z)t^2$ and $Y(z)$ respectively.)\\

\noindent
{\bf(e)} Special solutions

The determinant formulae of the special solutions are given as  
\begin{equation}\label{eq:E6solap}
\begin{array}{l}
\dfrac{1-\sfrac{f}{a_1}}{1-\sfrac{f}a_2}=\dfrac{a_1\prod_{i=1}^3(1-\sfrac{b_i}{a_1})}{a_2\prod_{i=1}^3(1-\sfrac{b_i}{a_2})}\dfrac{T_{a_1}(\tau_{m,n+1})T_{a_1}^{-1}(\tau_{m+1,n})}{T_{a_2}(\tau_{m,n+1})T_{a_2}^{-1}(\tau_{m+1,n})},
\\[5mm]
\dfrac{1-\sfrac{1}{b_2g}}{1-\sfrac{1}b_3g}=\dfrac{b_2\prod_{i=2}^3(1-\sfrac{a_i}{b_2})}{b_3\prod_{i=2}^3(1-\sfrac{a_i}{b_3})}\dfrac{T_{b_2}^{-1}(\tau_{m,n+1})T_{b_2}(\o{\tau}_{m+1,n})}{T_{b_3}^{-1}(\tau_{m,n+1})T_{b_3}(\o{\tau}_{m+1,n})}.
\end{array}
\end{equation}
Here the element $p_k$ in the determinant $\tau_{m,n}$ (\ref{eq:notationsap}) is given by
\begin{equation}\label{eq:E6pap}
p_k=\dfrac{b_3^k \left(\dfrac{a_3}{b_3};q\right)_k}{(q;q)_k} \varphi_D^{(2)} \Big(q^{-k},\frac{a_1}{b_1},\frac{a_2}{b_2},q^{-k+1}\frac{b_{3}}{a_{3}};q\frac{b_1}{a_{3}},q\frac{b_2}{a_{3}}\Big)
\end{equation}
and $\varphi_D^{(l)}$ is the $q$-Appell Lauricella function (i.e. the multivariable $q$- hypergeometric function)  \cite{GaR04} defined by
\begin{equation}\label{eq:appell}
\begin{array}l
\varphi_D^{(l)}(\alpha,\beta_1,\ldots,\beta_l,\gamma;z_1,\ldots,z_l):=\displaystyle\sum_{{m_i} \geq 0}\dfrac{(\alpha)_{|m|}(\beta_1)_{m_1}\ldots(\beta_l)_{m_l}}{(\gamma)_{|m|}(q)_{m_1}\ldots(q)_{m_l}}z_1^{m_1}\ldots z_l^{m_l}
\end{array}
\end{equation} 
where $|m|=m_1+\ldots+m_{l}$.    
      
 \begin{rem}\label{rem:trans}
{\rm {\bf On the transformation between the terminating $\varphi_D^{(2)}$ and the terminating ${}_3\varphi_2$}

\noindent     
The element $p_k$ (\ref{eq:E6pap}) expressed by the terminating $q$-Appell Lauricella series $\varphi_D^{(2)}$ \label{eq:appell} is rewritten in terms of the terminating $q$-hypergeometric series ${}_3\varphi_2$ (\ref{eq:qHGFap}) (big $q$-Jacobi polynomials \cite{KMNOY04, KS98, Noumi07}) under the transformation of $b_3=q^{k-1} a_1 a_2$ and $a_3=b_1 b_2$. $\square$
}
\end{rem}
These determinant formulae of the $q$-hypergeometric solutions (\ref{eq:E6solap}) are expected to be equivalent to those in \cite{Ikawa13, Nagao15}. (Note that there is a typographical error in eq.(38) in \cite{Ikawa13}, namely the expression $T_{a_2}T_{a_3}(\tau_{m,n-1})$ should read $T_{a_2}T_{a_4}(\tau_{m,n-1})$.) Determinant formulae of $q$-hypergeometric solutions still have not been given in terms of the non-terminating $q$-hypergeometric series ${}_3\varphi_2$ as far as we know.


\subsection{Case $q$-$D_5^{(1)}$}\label{subsec:D5ap}　\\

\noindent
The contents of these subsections (a), (b), (c) are the same as \cite{Ikawa13}.\\

\noindent
{\bf(a)} Setting of the Pad\'e approximation problem\\

In Table (\ref{eq:Ylistap}) the generating function is established as 
\begin{equation}\label{eq:D5Yap}
Y(x):=\dfrac{(a_1 x, a_2 x;q)_\infty}{(b_1 x, b_2 x;q)_\infty}
\end{equation}
and in Table (\ref{eq:Tlistap}) the time evolution is chosen as
\begin{equation}\label{eq:D5Tap}
T: (a_1,a_2,b_1,b_2,m,n) \mapsto (qa_1,a_2,qb_1,b_2,m,n).
\end{equation}\\

\noindent
{\bf(b)} Contiguity relations\\

By the definition (\ref{eq:GKHap}) we have the basic quantities
\begin{equation}\label{eq:D5GKHap}
G(x)=\displaystyle\prod_{i=1}^{2}\dfrac{(1-b_ix)}{(1-a_ix)},\quad 
K(x)=\dfrac{1-b_1x}{1-a_1x},\quad H(x)=\prod_{i=1}^{2}(1-a_ix),
\end{equation}
and by the expression (\ref{eq:Drelationap}) we obtain the Casorati determinants
\begin{equation}\label{eq:D5Dap}
\begin{array}{l}
D_1(x)=:\dfrac{c_0(1-xf)x^{m+n+1}Y(x)}{\prod_{i=1}^{2}(1-a_ix)},\quad
D_2(x)=:\dfrac{c_1x^{m+n+1}Y(x)}{1-a_1x},\quad D_3(x)=:\dfrac{c_2(1-b_1 x)x^{m+n+1}Y(x)}{\prod_{i=1}^{2}(1-a_ix)}
\end{array}
\end{equation}
where $f, c_0, c_1$ and $c_2$ are constants depending on parameters $a_1, a_2, b_1, b_2 \in \C^{\times}$, $m,n\in\Z_{\geq 0}$ but independent of $x$. Then the contiguity relations $L_2=0$ and $L_3=0$ are expressed by
\begin{equation}\label{eq:D5L2L3ap}
\begin{array}{l}
L_2(x)=C_0(1-xf)\o{y}(x)-(1-a_2x)y(qx)+(1-b_1x)\sfrac{y(x)}{g},\\
L_3(x)=C_1(1-\sfrac{x\o{f}}{q})y(x)+\sfrac{(1-a_1x)\o{y}(x)}{g}-q^{m+n+1}(1-\sfrac{b_2x}{q})\o{y}(\sfrac{x}{q})
\end{array}
\end{equation}
where $C_0=\sfrac{c_0}{c_1}, C_1=\sfrac{\o{c}_0}{c_1}$ and $g=\sfrac{c_1}{c_2}$.\\

Take note that in the items (c) and (d) we study the contiguity relations $L_2=0$ and $L_3=0$ (\ref{eq:D5L2L3ap}) for generic complex parameters $a_0=q^m$, $b_0=q^n$ $(m, n \in \C^{\times})$ and generic variables $f, g$ (depending on parameters $a_i, b_i \in \C^{\times}, i=0,1,2$) apart from the Pad\'e approximation problem (\ref{eq:padeap}) with eqs.(\ref{eq:D5Yap}) and (\ref{eq:D5Tap}). (see Remark \ref{rem:genericap})\\

\noindent
{\bf(c)} The $q$-Painlev\'e equation\\

Compatibility of the contiguity relations $L_2=0$ and $L_3=0$ (\ref{eq:D5L2L3ap}) gives the evolution equations and the constraint on the product $C_0C_1$ as follows:
\begin{equation}\label{eq:D5eqap}
g\u{g}=\dfrac{1}{qa_0b_0}\dfrac{(f-a_1)(f-b_1)}{(f-a_2)(f-b_2)},\quad
f\o{f}=a_2b_2\dfrac{(g-\sfrac{b_1}{a_0a_2})(g-\sfrac{a_1}{b_0b_2})}{(g-1)(g-\sfrac{1}{qa_0b_0})}
\end{equation}
and
\begin{equation}\label{eq:D5C0C1ap}
C_0C_1=\sfrac{(1-g)(1-qa_0b_0g)}{g^2}.
\end{equation}
The evolution equations (\ref{eq:D5eqap}) are equivalent to the $q$-Painlev\'e equation of type $D_5^{(1)}$ given in \cite{Ikawa13, JS96, KMNOY04, KNY15, Nagao15}. The 8 singular points in coordinates $(f,g)$ are on the four lines $f=0$, $f=\infty$, $g=0$ and $g=\infty$ as follows:
\begin{equation}
\begin{array}{l}
(f,g)=(a_1,0), (b_1,0), (0,\sfrac{b_1}{a_0a_2}), (0,\sfrac{a_1}{b_0b_2}),
(\infty,1), (\infty,\sfrac{1}{qa_0b_0}), (a_2,\infty), (b_2,\infty).
\end{array}
\end{equation}

\noindent
{\bf(d)} The Lax pair\\

The contiguity relations $L_2=0$ and $L_3=0$ (\ref{eq:D5L2L3ap}) give two scalar Lax equations $L_1=0$ and $L_2=0$ expressed by
\begin{equation}\label{eq:D5L1L2ap}
\begin{array}{l}
L_1(x)=\Big[(1-g) \left(1-qa_0b_0g\right)-\dfrac{x
   \left(a_1-b_0b_2 g\right)
   \left(b_1-a_0a_2 g\right)}{f}\Big]y(x)
   \\[5mm]
   \phantom{L_1(x): }
   +\dfrac{g
   (1-a_1 x) (1-a_2 x)}{1-f x}\Big[y(q
   x)-\dfrac{1-b_1 x}{g
   (1-a_2 x)} y(x)\Big]
   \\[5mm]
  \phantom{L_1(x): }
   +\dfrac{qa_0b_0g
   (1-b_1 x/q) (1-b_2 x/q)}{1-f x/q}
   \Big[y\left(\sfrac{x}{q}\right)-\dfrac{g
   (1-a_2 x/q)}{1-b_1
   x/q}y(x)\Big],
   \\[5mm]
L_2(x)=(1-xf)\o{y}(x)-(1-a_2x)y(qx)+\sfrac{(1-b_1x)y(x)}{g}.
\end{array}
\end{equation}
The scalar Lax pair $L_1=0$ and $L_2=0$ (\ref{eq:D5L1L2ap}) is equivalent to the 2 $\times$ 2 matrix ones in \cite{JS96, OWF11} and the scalar one in \cite{Nagao15, Yamada11} by using suitable gauge transformations of $y(x)$.\\

\noindent
{\bf(e)} Special solutions

The determinant formulae of the special solutions are given as
\begin{equation}\label{eq:D5solap}
\begin{array}{l}
\dfrac{1-\sfrac{f}{a_1}}{1-\sfrac{f}a_2}=\dfrac{a_1\prod_{i=1}^2(1-\sfrac{b_i}{a_1})}{a_2\prod_{i=1}^2(1-\sfrac{b_i}{a_2})}\dfrac{T_{a_1}(\tau_{m,n+1})T_{a_1}^{-1}(\tau_{m+1,n})}{T_{a_2}(\tau_{m,n+1})T_{a_2}^{-1}(\tau_{m+1,n})},
\\[5mm]
g=\dfrac{a_1(1-\sfrac{b_1}{a_1})}{q^n a_2 (1-\sfrac{b_2}{a_2})}\dfrac{T_{a_1}(\tau_{m,n+1})T_{a_1}^{-1}(\o{\tau}_{m+1,n})}{T_{a_2}(\o{\tau}_{m,n+1})T_{a_2}^{-1}(\tau_{m+1,n})}
\end{array}
\end{equation}
where the element $p_k$ in the determinant $\tau_
{m,n}$ (\ref{eq:notationsap}) is given by
\begin{equation}\label{eq:D5pap}
p_{k}=b^k_2\dfrac{(\sfrac{a_2}{b_2};q)_k}{(q;q)_k}{}_2\varphi_1 \Big(\substack{\displaystyle{q^{-k},\sfrac{a_1}{b_1}}\\[5mm]{\displaystyle{\sfrac{b_2q^{-k+1}}{a_2}}}};q,\dfrac{b_1}{a_2}q\Big).
\end{equation} 
The element $p_k$ (\ref{eq:D5pap}) is expressed in terms of the terminating $q$-hypergeometric series ${}_2\varphi_1$ (\ref{eq:qHGFap}) (little $q$-Jacobi polynomials \cite{KMNOY04, KS98, Noumi07}). These determinant formulae of the $q$-hypergeometric solutions (\ref{eq:D5solap}) are expected to be equivalent to those in \cite{Nagao15} and the terminating case of those in \cite{Sakai98}.\\

\subsection{Case $q$-$A_4^{(1)}$}
\label{subsec:A4ap}　\\

\noindent
{\bf(a)} Setting of the Pad\'e approximation problem\\

In Table (\ref{eq:Ylistap}) the generating function is established as 
\begin{equation}\label{eq:A4Yap}
Y(x):=\dfrac{(a_1x,a_2x;q)_\infty}{(b_1 x;q)_\infty}
\end{equation}
and in Table (\ref{eq:Tlistap}) the time evolution is chosen as
\begin{equation}\label{eq:A4Tap}
T: (a_1,a_2,b_1,m,n) \mapsto (qa_1,a_2,qb_1,m,n).
\end{equation}\\

\noindent
{\bf(b)} Contiguity relations

By the definition (\ref{eq:GKHap}) we have the basic quantities
\begin{equation}\label{eq:A4GKHap}
G(x)=\dfrac{1-b_1x}{(1-a_1x)(1-a_2x)},\quad K(x)=\dfrac{1-b_1x}{1-a_1x},\quad H(x)=(1-a_1x)(1-a_2x),
\end{equation}
and by the expression (\ref{eq:Drelationap}) we obtain the Casorati determinants
\begin{equation}\label{eq:A4Dap}
\begin{array}{l}
D_1(x)=:\dfrac{c_0(1-xf)x^{m+n+1}Y(x)}{(1-a_1x)(1-a_2x)},\quad
D_2(x)=:\dfrac{c_1x^{m+n+1}Y(x)}{1-a_1x},\quad D_3(x)=:\dfrac{c_2x^{m+n+1}Y(x)}{(1-a_1x)(1-a_2x)}
\end{array}
\end{equation}
where $f, c_0, c_1$ and $c_2$ are constants depending on parameters $a_1,a_2,b_1\in\C^{\times}$, $m,n\in\Z_{\geq 0}$ but independent of $x$. Then the contiguity relations $L_2=0$ and $L_3=0$ are expressed by
\begin{equation}\label{eq:A4L2L3ap}
\begin{array}{l}
L_2(x)=C_0(1-xf)\o{y}(x)-(1-a_2x)y(qx)+\sfrac{(1-b_1x)y(x)}{g},\\
L_3(x)=C_1(1-x\o{f}/q)y(x)+\sfrac{(1-a_1x)\o{y}(x)}{g}-q^{m+n+1}\o{y}(x/q)
\end{array}
\end{equation}
where $C_0=\sfrac{c_0}{c_1}, C_1=\sfrac{\o{c}_0}{c_1}$ and $g=\sfrac{c_1}{c_2}$.\\

Take note that in the items (c) and (d) below we study the contiguity relations $L_2=0$ and $L_3=0$ (\ref{eq:A4L2L3ap}) for generic complex parameters $a_0=q^m$, $b_0=q^n$ $(m, n \in \C^{\times})$ and generic variables $f, g$ (depending on parameters $a_0, a_1, a_2, b_0, b_1 \in \C^{\times}$) apart from the Pad\'e approximation problem (\ref{eq:padeap}) with eqs.(\ref{eq:A4Yap}) and (\ref{eq:A4Tap}). (see Remark \ref{rem:genericap})\\

\noindent
{\bf(c)} The $q$-Painlev\'e equation\\

Compatibility of the contiguity relations $L_2=0$ and $L_3=0$  (\ref{eq:A4L2L3ap}) gives the evolution equations and the constraint on the product $C_0C_1$ as follows:
\begin{equation}\label{eq:A4eqap}
g\u{g}=\dfrac{1}{qa_0b_0}\dfrac{(f-a_1)(f-b_1)}{f(f-a_2)},\quad
f\o{f}=-\dfrac{a_1a_2}{b_0}\dfrac{g-\sfrac{b_1}{a_0a_2}}{(g-1)(g-\sfrac{1}{qa_0b_0})}.
\end{equation}
and
\begin{equation}\label{eq:A4C0C1ap}
C_0C_1=\sfrac{(1-g)(1-qa_0b_0g)}{g^2}.
\end{equation}
The evolution equations (\ref{eq:A4eqap}) are equivalent to the $q$-Painlev\'e equation of type $A_4^{(1)}$ given in \cite{KMNOY04, KNY15, KTGR00, Nagao15}. The 8 singular points in coordinates $(f,g)$ are on the four lines $f=0, f=\infty, g=0$ and $g=\infty$ as follows:
\begin{equation}
\begin{array}{l}
(f,g)=(a_1,0), (b_1,0), (0,\sfrac{b_1}{a_0a_2}), (\infty,1), (\infty,\sfrac{1}{qa_0b_0}), (a_2,\infty),(\epsilon,-\sfrac{a_1}{b_0\epsilon})_2
\end{array}
\end{equation}
where the last point is a double point at $(0,\infty)$ with the gradient $fg=-\sfrac{a_1}{b_0}$. (The meaning of the double point is also written in \cite{KNY15}.)\\

\noindent
{\bf(d)} The Lax pair\\

The contiguity relations $L_2=0$ and $L_3=0$ (\ref{eq:A4L2L3ap}) give two scalar Lax equations $L_1=0$ and $L_2=0$ expressed by
\begin{equation}\label{eq:A4L1L2ap}
\begin{array}{l}
L_1(x)=\dfrac{qa_0b_0g(1 - b_1 x/q)}{f x/q-1}\left[y(\sfrac{x}{q})-\dfrac{g(1-a_2x/q)}{1-b_1x/q}y(x)\right]
 \\
\phantom{L_1(x):}+\dfrac{g(1-a_1x)(1-a_2x)}{f x-1}\left[y(q x)-\dfrac{1-b_1x}{g(1-a_2x)}y(x)\right]
 \\
\phantom{L_1(x):}+\left[\dfrac{a_1(b_1-a_0a_2g)x}{f} -(g-1)(qa_0b_0-1)\right]y(x),
\\
L_2(x)=(1-xf)\o{y}(x)-(1-a_2x)y(qx)+\sfrac{(1-b_1x)y(x)}{g}.
\end{array}
\end{equation}
The scalar Lax pair $L_1=0$ and $L_2=0$ (\ref{eq:A4L1L2ap}) is equivalent to the 2 $\times$ 2 matrix one for the $q$-Painlev\'e equation of type $q$-$P(A_4)$ in \cite{Murata09} and the scalar one in \cite{Nagao15} by using a suitable gauge transformation of $y(x)$.\\

\noindent
{\bf(e)} Special solutions

The determinant formulae of the special solutions are given as 
\begin{equation}\label{eq:A4solap}
\begin{array}{l}
\dfrac{1-\sfrac{f}{a_1}}{1-\sfrac{f}a_2}=\dfrac{a_1(1-\sfrac{b_1}{a_1})}{a_2(1-\sfrac{b_1}{a_2})}\dfrac{T_{a_1}(\tau_{m,n+1})T_{a_1}^{-1}(\tau_{m+1,n})}{T_{a_2}(\tau_{m,n+1})T_{a_2}^{-1}(\tau_{m+1,n})},
\\
g=\dfrac{a_1}{q^n a_2}\dfrac{T_{a_1}(\tau_{m,n+1})T_{a_1}^{-1}(\o{\tau}_{m+1,n})}{T_{a_2}(\o{\tau}_{m,n+1})T_{a_2}^{-1}(\tau_{m+1,n})}
\end{array}
\end{equation}
where the element $p_k$ in the determinant $\tau_
{m,n}$ (\ref{eq:notationsap}) is given by
\begin{equation}\label{eq:A4pap}
p_{k}=q^{\left(\substack{k\\2}\right)}\dfrac{(-a_2)^k}{(q;q)_k}{}_2\varphi_1 \Big(\substack{\displaystyle{q^{-k},\sfrac{a_1}{b_1}}\\[5mm]{\displaystyle{0}}};q,\dfrac{b_1}{a_2}q\Big).
\end{equation}
The element $p_k$ (\ref{eq:A4pap}) is expressed in terms of the terminating $q$-hypergeometric series ${}_2\varphi_1$ (\ref{eq:qHGFap}) ($q$-Laguerre polynomials
 \cite{KMNOY04, KS98, Noumi07}). These determinant formulae of the $q$-hypergeometric solutions (\ref{eq:A4solap}) are expected to be equivalent to those in \cite{Nagao15} and the terminating case of those in \cite{HK07}.\\
 
\subsection{Case $q$-$(A_2+A_1)^{(1)}$}\label{subsec:A21ap}　\\

\noindent
{\bf(a)} Setting of the Pad\'e approximation problem\\

In Table (\ref{eq:Ylistap}) the generating function is established as
\begin{equation}\label{eq:A21Yap}
Y(x):=(a_1x, a_2x;q)_{\infty}
\end{equation}
and in Table (\ref{eq:Tlistap}) the time evolution is chosen as
\begin{equation}\label{eq:A21Tap}
T: (a_1,a_2,m,n) \mapsto (qa_1,a_2,m,n).
\end{equation}\\

\noindent
{\bf(b)} Contiguity relations

By the definition (\ref{eq:GKHap}) we have the basic quantities
\begin{equation}\label{eq:A21GKHap}
G(x)=\dfrac{1}{(1-a_1x)(1-a_2x)},\quad K(x)=\dfrac{1}{1-a_1x},\quad H(x)=(1-a_1x)(1-a_2x),
\end{equation}
and by the expression (\ref{eq:Drelationap}) we obtain the Casorati determinants
\begin{equation}\label{eq:A21Dap}
\begin{array}{l}
D_1(x)=:\dfrac{c_0(1-xf)x^{m+n+1}Y(x)}{(1-a_1x)(1-a_2x)},\quad
D_2(x)=:\dfrac{c_1x^{m+n+1}Y(x)}{1-a_1x},\quad D_3(x)=:\dfrac{c_2x^{m+n+1}Y(x)}{(1-a_1x)(1-a_2x)}
\end{array}
\end{equation}
where $f, c_0, c_1$ and $c_2$ are constants depending on parameters $a_1, a_2\in\C^{\times}$, $m,n\in\Z_{\geq 0}$ but independent of $x$. Then the contiguity relations $L_2=0$ and $L_3=0$ are expressed by
\begin{equation}\label{eq:A21L2L3ap}
\begin{array}{l}
L_2(x)=C_0(1-xf)\o{y}(x)-(1-a_2x)y(qx)+\sfrac{y(x)}{g},\\
L_3(x)=C_1(1-x\o{f}/q)y(x)+\sfrac{(1-a_1x)\o{y}(x)}{g}-q^{m+n+1}\o{y}(x/q)
\end{array}
\end{equation}
where $C_0=\sfrac{c_0}{c_1}, C_1=\sfrac{\o{c}_0}{c_1}$and $g=\sfrac{c_1}{c_2}$.\\

Take note that in the items (c) and (d) below we study the contiguity relations $L_2=0$ and $L_3=0$ (\ref{eq:A21L2L3ap}) for generic complex parameters $a_0=q^m$, $b_0=q^n$ $(m, n \in \C^{\times})$ and generic variables $f, g$ (depending on parameters $a_0, a_1, a_2, b_0 \in \C^{\times}$) apart from the Pad\'e approximation problem (\ref{eq:padeap}) with eqs.(\ref{eq:A21Yap}) and (\ref{eq:A21Tap}). (see Remark \ref{rem:genericap})\\

\noindent
{\bf(c)} The $q$-Painlev\'e equation\\

Compatibility of the contiguity relations $L_2=0$ and $L_3=0$ (\ref{eq:A21L2L3ap}) gives the evolution equations and the constraint on the product $C_0C_1$ as follows:
\begin{equation}\label{eq:A21eqap}
g\u{g}=\dfrac{1}{qa_0b_0}\dfrac{f-a_1}{f-a_2},\quad
f\o{f}=-\dfrac{a_1a_2}{b_0}\dfrac{g}{(g-1)(g-\sfrac{1}{qa_0b_0})}
\end{equation}
and
\begin{equation}\label{eq:A21C0C1ap}
C_0C_1=\sfrac{(1-g)(1-qa_0b_0g)}{g^2}.
\end{equation}
The evolution equations (\ref{eq:A21eqap}) are equivalent to the $q$-Painlev\'e equation of type $(A_2+A_1)^{(1)}$, namely $P_{\rm IV}$, given in \cite{KNY01, KNY15, Nagao15, Sakai01}. The 8 singular points in coordinates $(f,g)$ are on the four lines $f=0, f=\infty, g=0$ and $g=\infty$ as follows:
\begin{equation}
\begin{array}{l}
(f,g)=(a_1,0), (\infty,1), (\infty,\sfrac{1}{qa_0b_0}), (a_2,\infty), (\epsilon,-\sfrac{a_1}{b_0\epsilon})_2, (\epsilon, -\sfrac{\epsilon}{a_0a_2})_2.
\end{array}
\end{equation}
Here the fifth point is a double point at $(0, \infty)$ with the gradient $fg=-\sfrac{a_1}{b_0}$ and the sixth point is a double point at $(0,0)$ with the gradient $\sfrac{g}{f}=-\sfrac{1}{a_0a_2}$. (The meaning of the two double points are also written in \cite{KNY15}.)\\

\noindent
{\bf(d)} The Lax pair\\

The contiguity relations $L_2=0$ and $L_3=0$ (\ref{eq:A21L2L3ap}) give two scalar Lax equations $L_1=0$ and $L_2=0$ expressed by
\begin{equation}\label{eq:A21L1L2ap}
\begin{array}{l}
L_1(x)=\dfrac{qa_0b_0g}{1-f x/q}\left[y(\sfrac{x}{q})-g(1-\sfrac{a_2x}{q})y(x)\right]
 \\
\phantom{L_1(x):}+\dfrac{g(1-a_1x)(1-a_2x)}{1-f x}\left[y(q x)-\dfrac{1}{g(1-a_2x)}y(x)\right]
 \\
\phantom{L_1(x):}+\left[\dfrac{a_0a_1a_2gx}{f} +(g-1)(qa_0b_0-1)\right]y(x),\\
L_2(x)=(1-xf)\o{y}(x)-(1-a_2x)y(qx)+\sfrac{y(x)}{g}.
\end{array}
\end{equation}
The scalar Lax pair $L_1=0$ and $L_2=0$ (\ref{eq:A21L1L2ap}) is equivalent to the 2 $\times$ 2 matrix one for the $q$-Painlev\'e equation of type $q$-$P(A_5)^{\#}$ in \cite{Murata09} and the scalar one in \cite{Nagao15} by using a suitable gauge transformation of $y(x)$.\\

\noindent
{\bf(e)} Special solutions

The determinant formulae of the special solutions are given as 
\begin{equation}\label{eq:A21solap}
\begin{array}{l}
\dfrac{1-\sfrac{f}{a_1}}{1-\sfrac{f}a_2}=\dfrac{a_1}{a_2}\dfrac{T_{a_1}(\tau_{m,n+1})T_{a_1}^{-1}(\tau_{m+1,n})}{T_{a_2}(\tau_{m,n+1})T_{a_2}^{-1}(\tau_{m+1,n})},
\\[5mm]
g=\dfrac{a_1}{q^n a_2}\dfrac{T_{a_1}(\tau_{m,n+1})T_{a_1}^{-1}(\o{\tau}_{m+1,n})}{T_{a_2}(\o{\tau}_{m,n+1})T_{a_2}^{-1}(\tau_{m+1,n})}
\end{array}
\end{equation}
where the element $p_k$ in the determinant $\tau_
{m,n}$ (\ref{eq:notationsap}) is given by
\begin{equation}\label{eq:A21pap}
p_{k}=(-1)^kq^{\left(\substack{k\\2}\right)}\dfrac{(a_2;q)_k}{(q;q)_k}{}_1\varphi_1 \Big(\substack{\displaystyle{q^{-k}}\\[5mm]{\displaystyle{0}}};q,\dfrac{a_1}{a_2}q\Big)
\end{equation}
The element $p_k$ (\ref{eq:A21pap}) is expressed in terms of the terminating $q$-hypergeometric series ${}_2\varphi_1$ (\ref{eq:qHGFap}) (Stieltjes-Wigert polynomials \cite{KMNOY04, KS98, Noumi07}). These determinant formulae of the $q$-hypergeometric solutions (\ref{eq:A21solap}) are expected to be  equivalent to those in \cite{Nagao15} and the terminating case of those in \cite{Nakazono10}.\\


\section{Conclusion}\label{sec:conc}
\subsection{Summary}　\\

\noindent
In this paper for the generating function $Y(x)$ given in Table (\ref{eq:Ylistap}), we established the Pad\'e approximation problem related to the $q$-Painlev\'e equations of type $E_6^{(1)}$, $D_5^{(1)}$, $A_4^{(1)}$ and $(A_2+A_1)^{(1)}$. Then for the time evolution $T$ given in Table (\ref{eq:Tlistap}), we established another Pad\'e approximation problem. By solving these problems, we derived the evolution equations, the scalar Lax pairs and the determinant formulae of the special solutions for the corresponding $q$-Painlev\'e equations. The main results are given in Section \ref{sec:approximation results}. 
\subsection{Problems}　\\

\noindent
Some open problems related to the results of this paper are as follows:

{\bf 1.} In this paper by choosing certain time evolutions $T$, we applied the Pad\'e method to each type $q$-$E_6^{(1)}$, $q$-$D_5^{(1)}$, $q$-$A_4^{(1)}$ and $q$-$(A_2+A_1)^{(1)}$. Moreover by choosing other time evolutions, we can perform similar computations. It may be interesting to study 
the relation between the Pad\'e approximation method for the various time evolutions and Ba\"cklund transformations of the affine Weyl group, for example the case $q$-$E_6^{(1)}$ of the interpolation problem in \cite{Ikawa13}.

{\bf 2.} By the results of this paper, it turned out that the Pad\'e approximation method could be applied to the $q$-Painlev\'e equations of type $E_6^{(1)}$, $D_5^{(1)}$, $A_4^{(1)}$ and $(A_2+A_1)^{(1)}$. It may be interesting to study the degenerations between these results.

{\bf 3.} In this paper we applied the Pad\'e method of differential grid (i.e. Pad\'e approximation) to the $q$-Painlev\'e equations of type $E_6^{(1)}$, $D_5^{(1)}$, $A_4^{(1)}$ and $(A_2+A_1)^{(1)}$.  It may be interesting to study whether the Pad\'e method can be also applied to additive difference Painlev\'e equations.

{\bf 4.} It may be interesting to study whether the Pad\'e method can be further applied to the other generalized Painlev\'e systems, for example the $q$-Garnier system in \cite{Sakai05} and the higher order $q$-Painlev\'e system in \cite{Suzuki15}.\\

\section*{Acknowledgment}

\noindent
The author is grateful to Professor Yasuhiko Yamada for valuable discussions on this research, and for his encouragement. The author also thanks Professors Nobutaka Nakazono, Masatoshi Noumi, Hidetaka Sakai, Takao Suzuki, Teruhisa Tsuda, and the referee for stimulating comments and for kindhearted support. 


\begin{thebibliography}{A}
\bibitem{Clarkson13}
Clarkson P.A., {\it Recurrence coefficients for discrete orthonormal polynomials and the Painlev\'e equations}, J. Phys., A {\bf 46}, no.18, (2013), 185205--185222.
\bibitem{GaR04}
Gasper G., and Rahman M., {\it Basic Hypergeometric Series. With a foreword by Richard Askey. Second edition. Encyclopedia of Mathematics and its Applications}, Volume {\bf 96}, Cambridge University Press, Cambridge, (2004). 
\bibitem{HK07}
Hamamoto T., and Kajiwara K., {\it Hypergeometric solutions to the q-Painlev\'e equation of type $A_4^{(1)}$}, J. Phys. A: Math. Theor., {\bf 40} (2007), 12509--12524. 
\bibitem{Ikawa13}
Ikawa Y., {\it Hypergeometric Solutions for the $q$-Painlev\'e Equation of Type $E_6^{(1)}$ by the Pad\'e method}, Lett. Math. Phys., Volume {\bf 103}, Issue 7  (2013), 743--763.
\bibitem{JS96}
Jimbo M., and Sakai H., {\it A $q$-analog of the sixth Painlev\'e equation}, Lett. Math. Phys., {\bf 38} (1996), 145--154.
\bibitem{KMNOY04}
Kajiwara K., Masuda T., Noumi M., Ohta Y., and Yamada Y.,
{\it Hypergeometric solutions to the $q$-Painlev\'e equations}, 
Int. Math. Res. Not. 2004, {\bf 47} (2004), 2497--2521.
\bibitem{KNY01}
Kajiwara K., Noumi M., and Yamada Y., {\it A study on the fourth $q$-Painlev\'e equation}, J. Phys. A: Math. Gen., {\bf 34} (2001), 8563--8581.
\bibitem{KNY15}
Kajiwara K., Noumi M., and Yamada Y.,
{\it Geometric aspects of Painlev\'e equations}, arXiv 1509.08186 [nlin.SI].
\bibitem{KS98}
Koekoek R., and Swarttouw R. F., {\it The Askey-scheme of hypergeometric orthogonal polynomials and its q-analogue}, Delft University of Technology, Department of Technical Mathematics and Informatics Report (1998) no. 98--17.
\bibitem{KTGR00}
Kruskal M. D., Tamizhmani K. M., Grammaticos B., and Ramani A., {\it Asymmetric discrete Painlev\'e equations}, Regul. Chaot. Dyn., {\bf 5} (2000), 273--281.
\bibitem{Magnus95}
Magnus A., {\it Painlev\'{e}-type differential equations for the recurrence coefficients of semi- classical orthogonal polynomials}, J. Comput. Appl. Math., {\bf 57} (1995), 215--237.
\bibitem{Murata09}
Murata M., {\it Lax forms of the $q$-Painlev\'e equations}, J. Phys. A: Math. Theor., {\bf 42} (2009), 115201--115217.
\bibitem{Nagao15}
Nagao, H., {\it The Pad\'e interpolation method applied to $q$-Painlev\'e equations}, Lett. Math. Phys. {\bf 105} (2015), no. 4, 503--521.
\bibitem{Nakazono10}
Nakazono N., {\it Hypergeometric $tau$ Functions of the $q$-Painlev\'e Systems of Type $(A_2+A_1)^{(1)}$} SIGMA, {\bf 6} (2010), 084--099.
\bibitem{Nakazono13}
Nakazono N., {\it Solutions to discrete Painlev\'e systems arising from two types of orthogonal polynomials}, Reports of RIAM Symposium, {\bf 23}AO--S7 (2013), 35--41.
\bibitem{Noumi07}
Noumi M., {\it Special functions arising from discrete Painlev\'e equations: a survey.}, J. Comput. Appl. Math. {\bf 202} (2007), no. 1, 48--55.
\bibitem{NTY13}
Noumi M., Tsujimoto S., and Yamada Y., {\it Pad\'e interpolation for elliptic Painlev\'e equation}, Symmetries, integrable systems and representations, Springer Proc. Math. Stat., Volume {\bf 40} (2013), 463--482.
\bibitem{OWF11}
Ormerod C.M., Witte N.S., and Forrester P.J., {\it Connection preserving deformations and q-semi-classical orthogonal polynomials}, Nonlinearity, {\bf 24} (2011), 2405--2434. 
\bibitem{RGTT01}
Ramani A., Grammaticos B., Tamizhmani T., and Tamizhmani K.M., Special Function Solutions of the Discrete Painlev\'e Equations, Comput. Math. Appl., {\bf42} (2001), no. 3--5, 603--614.
\bibitem{Sakai98}
Sakai H., {\it Casorati determinant solutions for the $q$-difference sixth Painlev\'e equations}, Nonlinearity, {\bf 11} (1998), 823--833.
\bibitem{Sakai01}
Sakai H., {\it Rational surfaces with affine root systems and geometry of the Painlev\'e equations}, Commun. Math. Phys., {\bf 220} (2001), 165--221.
\bibitem{Sakai05}
Sakai H., {\it A $q$-Analog of the Garnier System}, Funkcial. Ekvac. {\bf 48} (2005), 273--297.
\bibitem{Sakai06}
Sakai H., {\it Lax form of the $q$-Painlev\'e equation associated with the $A_2^{(1)}$ surface}, J. Phys. A: Math. Gen., {\bf 39} (2006), 12203--12210.
\bibitem{Suzuki15}
Suzuki T., {\it A $q$-analogue of the Drinfeld-Sokolov hierarchy of type $A$ and $q$-Painlev\'e system}, AMS Contemp. Math. {\bf 651} (2015), 25--38.
\bibitem{Tsuda10}
Tsuda T., {\it On an integrable system of q-difference equations satisfied by the universal characters: its Lax formalism and an application to $q$-Painlev\'e equations}, Comm. Math. Phys. {\bf 293} (2010), 347--359.
\bibitem{Van07}
Van W.A., {\it Discrete Painlev\'e equations for recurrence coefficients of orthogonal polynomials. Difference equations, special functions and orthogonal polynomials}, World Sci. Publ., Hackensack, NJ (2007), 687--725. 
\bibitem{Witte09}
Witte N.S., {\it Biorthogonal Systems on the Unit Circle, Regular Semiclassical Weights, and the Discrete Garnier Equations}, IMRN, {\bf 6}(2009), 988--1025. 
\bibitem{Witte15}
Witte N.S., {\it Semiclassical orthogonal polynomial systems on nonuniform lattices, deformations of the Askey table, and analogues of isomonodromy}, Nagoya Mathematical Journal, {\bf 219} (2015),127-234
\bibitem{WO12}
Witte N.S., and Ormerod C.M, {\it Construction of a Lax Pair for the $E_6^{(1)}$ $q$-Painlev\'e System}, SIGMA, {\bf 8} (2012), 097--123. 
\bibitem{Yamada09-1}
Yamada Y., {\it Pad\'e method to Painlev\'e equations}, Funkcial. Ekvac., {\bf 52} (2009), 83--92.
\bibitem{Yamada09-2} 
Yamada Y., {\it A Lax formalism for the elliptic difference Painlev\'e equation}, SIGMA {\bf 5} (2009) 042 (15pp).
\bibitem{Yamada11}
Yamada Y., {\it Lax formalism for $q$-Painlev\'e equations with affine Weyl group symmetry of type $E^{(1)}_n$}, IMRN, {\bf 17} (2011), 3823--3838.
\bibitem{Yamada14}
Yamada Y., {\it A simple expression for discrete Painlev\'e equations}, RIMS Kokyuroku Bessatsu, B{\bf47} (2014), 087--095.

\end{thebibliography}
\end{document}